\documentclass[12pt, reqno]{article}
\usepackage{amsmath, amssymb, amscd, eucal}
\usepackage{latexsym}
\usepackage{graphicx}

\newtheorem{thm}{Theorem}
\newtheorem{lem}{Lemma}
\newtheorem{defin}{Definition}
\newenvironment{prf}{\medskip 
\noindent {\em Proof.}}{\hfill $\square$ \\}
\newenvironment{rmk}{\medskip 
\noindent {\bf Remark.}}{\hfill \mbox{} \\}

\def\ms{\medskip}
\def\bs{\bigskip}

\def\nind{\noindent}

\def\eg{e.\ g.\ }
\def\ie{i.\ e.\ }

\def\ra{\rightarrow}

\def\({\left(}                   
\def\){\right)}                  
\def\[{\left[}                   
\def\]{\right]}                  
\def\lan{\langle}
\def\ran{\rangle}

\def\operatorname#1{{\rm #1\,}}                                
\def\op{\operatorname}

\def\tms{\times}                 
                  
\def\id{\equiv}                  
\def\grd{\nabla}         
\def\del{\partial}   
\def\lap{\triangle}

\def\a{\alpha}

\def\d{\delta}         
\def\e{\epsilon}         
            
\def\th{\theta}

\def\m{\mu}       
\def\n{\nu}

\def\s{\sigma}    
\def\ph{\phi}     
     
\def\om{\omega}

\def\Lam{\Lambda} 
       
\def\S{\Sigma}

\def\Om{\Omega}

\def\R{\mathbb R}

\def\H{\mathbb H}

\newcommand{\bg}{\overline g}
\newcommand{\bM}{\overline M}
\newcommand{\tg}{\tilde g}

\newcommand{\bB}{\overline B}
\newcommand{\bR}{\overline R}
\newcommand{\bH}{\overline H}

\title
{On the Uniqueness of the ADS Spacetime}
\author   {Xiaodong Wang 
\thanks{Department of Mathematics, MIT, Cambridge, MA 02139, USA. 
\tt{Email:xwang@math.mit.edu} }}
\date{September 5, 2002}

\begin{document}
\maketitle
\bibliographystyle{amsalpha}
\pagestyle{plain}
\begin{abstract}
The uniqueness of the ADS spacetime among all stactic vacuum spacetimes
with the same conformal infinity is proved in all dimensions under
the spin assumption. 
\end{abstract}

\nind

\nind
{\bf MR(2000) Subject Classification} 53, 58 

\bs
In Einstein's theory of general
relativity with a negative cosmological constant $\Lam$,
a vacuum spacetime is a solution to the equation $R_{ab}=\Lam g_{ab}$
and the lowest-energy solution is the anti-de Sitter spacetime. Moreover
it was proved by Bocher-Gibbons-Horowitz \cite{BGH} in $3+1$ dimensions
that the anti-de Sitter spacetime is the unique
static, asymptotically anti-de Sitter vacuum.

To give the precise statement of their theorem, let us first recall that
an $(n+1)$ dimensional static spacetime $(N, \bg)$ has the form
\begin{equation}\label{sta}
\begin{split}
N&=\R\tms M, \\
\bg&=-V^2dt^2+g,
\end{split}
\end{equation}
where $(M,g)$ is a Riemannian manifold and $V$ is a positive function
on $M$. The vacuum Einstein equation (without loss of generality
we always take the negative cosmological constant $\Lam$ to be $-n$) 
\begin{equation}\label{va}
\op{Ric}(\bg)=-n \bg
\end{equation}
can be written in terms of $g$ and $V$ as
\begin{align}
\op{Ric}(g)+ng&=V^{-1}D^2V, \label{f1} \\
\lap V&=nV,\label{f2}
\end{align}
where $D^2$ and $\lap$ are the Hessian and Laplacian on $(M,g)$. From
these two equations it is easy to see that $g$ has scalar curvature
$S=-n(n-1)$. We will often just call the triple $(M,g,V)$ a static vacuum.

The anti-de Sitter spacetime, in $(n+1)$ dimensions,  can be written
as 
$$ds^2=-(1+r^2)dt^2+(1+r^2)^{-1}dr^2+r^2d\om^2, $$ 
which is manifestly static. It is simply $\R\tms \H^n$
with the metric $-\cosh^2rdt^2+h$,
where $h=dr^2+\sinh^2rd\om^2$ is the hyperbolic metric in polar
coordinates. Another interesting example is the so called AdS soliton,
which is $\R\tms B^2\tms T^{n-2}$ with the metric
$$ds^2=-r^2dt^2+V(r)^{-1}dr^2+V(r)d\ph^2+r^2h,$$
where, $h$ is a flat metric on $T^{n-2}$,
$V(r)={r^2}\(1-\frac{r_0^n}{r^n}\)$, with 
$r_0>0$ a constant and $\ph$ is periodic with period $4\pi /nr_0$
to resolve the singularity at $r=r_0$. A uniqueness theorem has recently
been established by Galloway-Surya-Woolgar \cite{GSW}. In both
examples $(M, g)$ is conformally compact. The conformal infinity is
the standard sphere $S^{n-1}$ in the anti-de Sitter case and a flat
torus $T^n$ in the AdS soliton case.

The simplest example of a spacetime which is asymptotic to the anti-de Sitter
spacetime is the Schwarzschild-anti-de Sitter metric:
\begin{equation}\label{sads}
ds^2_0=-f(r)dt^2+f(r)^{-1}dr^2+r^2d\om^2,
\end{equation}
where $f(r)=1+r^2-M/r^{n-2}$. In dimension $n=3$ 
Boucher-Gibbons-Horowitz \cite{BGH} defines
a spacetime to be asymptotically anti-de Sitter if outside of a spatially 
compact world tube the metric has the following asymptotic behavior:
\begin{equation}\label{asp}
\begin{split}
ds^2&=ds^2_0+O(r^{-2})dt^2+O(r^{-6})dr^2
+O(r)(\text{ remaining differentials not involving } dr) \\
&+O(r^{-1})(\text{ remaining differentials involving } dr),
\end{split}
\end{equation}
where $ds^2_0$ is metric in (\ref{sads}). With this definition they
proved that the only static asymptotically anti-de Sitter solution to
$R_{ab}=\Lam g_{ab}$ is the anti-de Sitter spacetime.

There are two important ingredients in the
proof. The first one is the positive mass theorem for
asymptotically hyperbolic manifolds. The second ingredient is the
following remarkable identity attributed to Lindblom 
\begin{equation}\label{lin}
\op{div}\(V^{-1}\grd(W-W_0)\)=\frac{1}{4}V^3W^{-1}|B|^2
+\frac{3}{4}V^{-1}W^{-1}|\grd (W-W_0)|^2,
\end{equation}
where $W=|\grd V|^2$, $W_0=V^2-1$ and $B$ is the Bach tensor. In a local
frame the Bach tensor $B_{ijk}=R_{ij,k}-R_{ik,j}+\frac{1}{4}S_jh_{ik}
-\frac{1}{4}S_kh_{ij}$ and its importance lies in the fact that a
Riemannian $3-$manifold is locally conformally flat iff its Bach tensor $B=0$.
To verify the identity (\ref{lin}), 
the key point is
to express $B$ in terms of $V$ and $h$  and this is possible because
the curvature tensor is determined by the Ricci tensor
in dimension three. To see this
choose a local orthonormal frame and compute, using (\ref{f1}) and $S=-6$
\begin{align*}
B_{ijk}&=R_{ij,k}-R_{ik,j} \\
&=V^{-1}(V_{ij,k}-V_{ik,j}+V^{-2}\(V_jV_{ik}-V_kV_{ij}\) \\
&=V^{-1}R_{jkli}V_l+V^{-2}\(V_jV_{ik}-V_kV_{ij}\).
\end{align*}
In dimension three, the curvature tensor
$R_{jkli}=3(\d_{ji}\d_{kl}-\d_{jl}\d_{ki})+R_{jl}\d_{ki}+\d_{kl}R_{ki}
-R_{ji}\d_{kl}-\d_{ji}R_{kl}$. By using (\ref{f1}) again we end up
with
$$B_{ijk}=V^{-2}\[(V_jlV_l-3VV_j)\d_{ik}-(V_klV_l-3VV_k)\d_{ij}
+2(V_{ik}V_j-V_{ij}V_k)\].$$
In higher dimensions the curvature tensor is considerably more complicated
with an extra piece, the Weyl
tensor which is out of the reach of the field equations (\ref{f1})(\ref{f2}),
therefore it is not clear at all how to generalize (\ref{lin}).

The main purpose of this paper is to give a different approach which is
very elementary and works in any dimensions. We
also use the positive mass theorem, but in place of
the Lindblom identity (\ref{lin}) we use an elementary identity in
Riemannian geometry which holds in all dimensions.
Let $f$ be a function and $T$ a symmetric
$2-$tensor on a Riemannian manifold $(M,g)$, then we have
\begin{equation}
\lan T, D^2f\ran_g=\op{div}(i_{\grd f}T)-\lan df,\d T\ran_g.
\end{equation}
In a local orthonormal frame this can be written as
\begin{equation}
T_{ij}f_{ij}=\(T_{ij}f_i\)_j-f_iT_{ij,j}
\end{equation}
If $\Om\subset M$ is a compact and smooth domain, the by Stokes
theorem we obtain
\begin{equation}\label{int}
\int_{\Om}\lan T, D^2f\ran_gd\m=\int_{\del \Om}T(\grd f,\n)dA
-\int_{\Om}\lan df,\d T\ran_gd\m 
\end{equation}
where $\n$ is the outer unit normal of $\del\Om$.

An advantage of our approach is that we use conformal compactification
to formulate the asymptotics (see Wang \cite{Wang2} and 
Chrusiel-Herzlich \cite{CH}), which allows more general asymptotics. 
Instead of having mass $M$ as just a number, 
there is a {\em mass aspect} $\tau$ which is
a $2$-tensor. 
Moreover the asymptotics required for the application of the positive
mass theorem follow naturally from
the field equations (\ref{f1}) (\ref{f2}). With Fefferman and Graham's
result on conformally compact Einstein metrics the analysis is also 
made clean and transparent.
We state our main result as

\begin{thm}
Let $(M,g,V)$ be a static solution to the vacuum Einstein equation with
negative cosmological constant (\ref{va}). If
\begin{enumerate}
\item $(M,g)$ is conformally compact and the conformal boundary is $S^{n-1}$,
\item $V^{-1}$ is a defining function and $V^{-2}g|_{S^{n-1}}$ is the
      standard metric on $S^{n-1}$
\item $M$ is spin,
\end{enumerate}
then $(M,g,V)$ is the anti-de Sitter spacetime.
\end{thm}

\begin{rmk}
The last assumption that $M$ is spin, which is superfluous in dimension
$n=3$, is needed because that is the only case the positive mass theorem
for asymptotically hyperbolic manifolds is known. 
\end{rmk}

Since $(M, g)$ is conformally compact with scalar curvature $S=-n(n-1)$
we can choose a defining function $r$ such that near infinity
\begin{equation}\label{good}
g=r^{-2}(dr^2+h_r),
\end{equation}
where $h_r$ is an $r-$dependent family of metrics on the conformal boundary
$S^{n-1}$ with $h_0$ being the unit round metric (see \eg Lemma 2.1 in
Graham \cite{Gra}). We choose local coordinates $x^1,\ldots,x^{n-1}$ 
on the boundary and write $h_r=h_{ij}(r,x)dx^idx^j$.
Moreover in view of the assumption (ii) we can
assume $Vr|_{S^{n-1}}\id 1$. 

For example it is easy to see that
the anti-de Sitter spacetime can be written in this form as
\begin{align}
g_0&=r^{-2}\(dr^2+(1-r^2/4)^2h_0\) \label{ex1} \\
V&=1/r+r/4.\label{ex2}
\end{align}

\ms

We need to study the asymptotic expansion of $h_r$ and $V$ in detail.
To do this one could study directly 
the equations (\ref{f1}) and (\ref{f2}). Instead we take
advantage of the various results available for conformally compact Einstein
manifolds. We recommend
Graham \cite{Gra} for a clear and succinct introduction.

Consider the manifold $N=S^1\tms M$ with the metric
$$\tg=V^2d\th^2+g=r^{-2}\(dr^2+(Vr)^2d\th^2+h_r\),$$
where $\th$ is periodic with period $2\pi$. By the equations (\ref{f1})
(\ref{f2})
$(N,\tg)$ is conformally compact Einstein 
(\ie $\op{Ric}(\tg)=-n\tg$) with conformal
infinity $(S^1\tms S^{n-1}, d\th^2+h_0)$.
If $(M,g, V)$ is the anti-de Sitter spacetime, we have by (\ref{ex1})
(\ref{ex2})
$$\tg_0=r^{-2}\(dr^2+(1+r^2/4)^2d\th^2+(1-r^2/4)^2h_0\). $$

Both $\tg$ and $\tg_0$ are conformally compact Einstein with the same
conformal infinity $(S^1\tms S^{n-1}, d\th^2+h_0)$. By the work of
Fefferman and Graham as presented in \cite{Gra} (roughly speaking, in the
Taylor expansion in $r$ the first $n$ terms are locally determined by
the metric $h_0$ on the conformal boundary), we have
\begin{equation}\label{exp}
(Vr)^2d\th^2+h_r=(1+r^2/4)^2d\th^2+(1-r^2/4)^2h_0+r^n(\a d\th^2+\tau)
+o(r^{n}),
\end{equation}
where $\tau$ is a symmetric 2-tensor on $S^{n-1}$ and 
\begin{equation}\label{alp}
\a=-\op{tr}_{h_0}\tau.
\end{equation}

The equation (\ref{exp}) can be separated as two equations
\begin{align}
V&=1/r+r/4+\a r^{n-1}/2+o(r^{n-1}) \label{expv}\\
h_r&=(1-r^2/4)^2h_0+\tau r^n+o(r^n).\label{exph}
\end{align}
These asymptotic expansions can be differentiated. For example by
differentiating (\ref{expv}) in $r$ we obtain
\begin{equation}\label{dvr}
\frac{\del V}{\del r}=
-1/r^2+1/4+(n-1)\a r^{n-2}/2+o(r^{n-2}).
\end{equation}

The equation (\ref{exph}) shows that $(M, g)$ is asymptotically
hyperbolic in the sense of \cite{Wang2}. 

\begin{defin}\cite{Wang2}\label{AH}
A conformally compact  manifold $(X, g)$ is called asymptotically
hyperbolic if it satisfies:
\begin{enumerate}
\item the conformal infinity is the standard sphere $(S^{n-1}, h_0)$,
\item there exists a good defining function $r$ such that  
\begin{equation}
g=r^{-2}(dr^2+h_r)
\end{equation}
in a collar neighborhood of the conformal infinity and 
\begin{equation}
h_r=(1-r^2/4)^2h_0+\tau r^n+o(r^{n}),
\end{equation}
where $\tau$ is a symmetric 2-tenor on $S^{n-1}$. Moreover the asymptotic 
expansion can be differentiated twice.
\end{enumerate}
\end{defin}

For asymptotically hyperbolic spin manifolds, the following positive mass
theorem was proved in Wang \cite{Wang2} (see also Chrusciel-Herzlich 
\cite{CH}). 
\begin{thm}\cite{Wang2}\label{pm}
Let $(X,g)$ be an asymptotically hyperbolic manifold.
If $X$ is spin and has scalar curvature $R\geq -n(n-1)$
then we have 
$$\int_{S^{n-1}}\op{tr}_{h_0}(\tau)d\m_{h_0}\geq 
\left|\int_{S^{n-1}}\op{tr}_{h_0}(\tau)x d\m_{h_0}\right|.$$
Moreover equality holds iff $(X,g)$ is isometric to the hyperbolic
space $\H^n$.
\end{thm}

We use (\ref{expv})(\ref{exph})
and (\ref{good}) to compute 
\begin{align*}
|\grd V|^2-V^2+1&=r^2\[\(\frac{\del V}{\del r}\)^2
+h^{ij}\frac{\del V}{\del x^i}\frac{\del V}{\del x^j}\]-V^2+1 \\
&=r^2\[\(-1/r^2+1/4+(n-1)\a r^{n-2}/2+o(r^{n-2})\)^2+O(r^{2(n-1)})\] \\
&\quad -\(1/r+r/4+\a r^{n-1}/2+o(r^{n-1})\)^2+1 \\
&=-n\a r^{n-2}+o(r^{n-2}).
\end{align*}
To summarize we have shown
\begin{equation}\label{bm}
|\grd V|^2-V^2=-1-n\a r^{n-2}+o(r^{n-2})
\end{equation}
Differentiating in $r$ gives
\begin{equation}\label{dint}
\frac{\del}{\del r}\(|\grd V|^2-V^2\)=-n(n-2)\a r^{n-3}+o(r^{n-3})
\end{equation}

\nind
{\em Proof of the Theorem:}
Define $M^{\e}=\{p\in M|r(p)\geq \e\}$. For $\e$ small, this is a compact
and smooth domain in $M$. 
Let $T=\op{Ric}(g)+(n-1)g$ which is actually the Einstein tensor of $g$. 
By (\ref{f2})
we have
\begin{equation}\label{ein}
T=V^{-1}D^2V-g.
\end{equation}
Hence $V|T|^2=\lan D^2V-Vg, T\ran_g=\lan D^2V,T\ran_g$.
By the second Bianchi identity $\d T=(\frac{1}{2}-\frac{1}{n})dS=0$ 
for $S=-n(n-1)$.
Therefore by the formula (\ref{int})
\begin{equation}
\int_{M^{\e}}V|T|_g^2d\m_g=\int_{\del M^{\e}}T(\grd V, \n)d\s.
\end{equation}
The outer unit normal $\n=-r\frac{\del}{\del r}$ and the induced area form
$d\s=r^{-(n-1)}\sqrt{H(r,x)}dx$ with $H(r,x)=\op{det}(h_{ij}(r,x))$. Therefore,
using (\ref{ein}) again, we can rewrite the above equation as
\begin{align*}
\int_{M^{\e}}V|T|_g^2d\m_g
&=-\e^{-(n-2)}\int_{S^{n-1}}\(V^{-1}D^2V(\grd V,\frac{\del}{\del r})
-\lan \grd V,\frac{\del}{\del r}\ran\) \sqrt{H(\e,x)}dx \\
&=-\frac{1}{2}\e^{-(n-2)}\int_{S^{n-1}}V^{-1}\frac{\del}{\del
r}\(|\grd V|^2-V^2\) \sqrt{H(\e,x)}dx \\
&=\frac{n(n-2)}{2}\int_{S^{n-1}}\a d\s_{h_0}+o(1),
\end{align*}
where in the last step we use the asymptotic expansions (\ref{expv}) and
(\ref{dint}). Taking
$\e\ra 0$, in view of (\ref{alp}),  we obtain
$$\int_{M}V|T|_g^2d\m_g=-\frac{n(n-2)}{2}\int_{S^{n-1}}
\op{tr}_{h_0}(\tau) dA_{h_0}.$$
But by Theorem \ref{pm}, the right hand side is non-positive.
Therefore both sides must be zero. By the characterization of zero
mass case in Theorem \ref{pm}, 
$(M,g)$ is the hyperbolic space $\H^n$. Then it is an
easy matter to show that the triple $(M,g,V)$ is the anti-de
Sitter spacetime.
\hfill $\square$ \\

In the last step we only need the positive mass theorem to the extent
that the mass is nonnegative and the characterization of zero mass case
can be replaced by an elementary argument.
In fact we have $T=0$, \ie $\op{Ric}(g)+(n-1)g=V^{-1}(D^2V-Vg)=0$. This
implies that $\bg=V^{-2}g$ is Einstein. Therefore $(\bM, \bg)$ is a
compact Einstein manifold with a totally geodesic boundary which is the
standard sphere $S^{n-1}$. Moreover $\overline{\grd}V$ is a conformal 
vector field. It is then elementary to show that $(\bM, \bg)$ is the 
standard hemisphere $S^n_+\subset R^{n+1}$ and $V(x)=1/x_{n+1}$. This 
easily implies that $(M,g,V)$ is the anti-de Sitter spacetime.

\ms

{\bf Acknowledgment: } I am indebted to Professor Galloway for stimulating
conversations, from whom I learned this problem during the AIM-Stanford
workshop on general relativity in May, 2002. I thank the workshop
organizers, especially Professor Rick Schoen, for inviting me and for
providing an excellent environment. I also wish to thank Professor Tian
for constant encouragement.
\bs
\begin{center}
------------------------------------------------------------
\end{center}
\ms

After the paper was submitted, the paper Qing \cite{Qing} appeared.
Inspired by his idea, I give
an alternative approach which does not require the
spin assumption in dimension $n\leq 7$. It is interesting to note
that this new approach  uses the positive mass theorem for
asymptotically flat manifolds while the previous approach uses
the positive mass theorem for asymptotically hyperbolic manifolds.

Let $(M,g, V)$ be a static AdS vacuum spacetime.
By the equations (\ref{f1}) (\ref{f2}) and the Bochner formula we compute
\begin{align*}
\ \frac{1}{2}\lap\(|\grd V|^2-V^2\)
&=|D^2V|^2+\grd V\cdot\grd\lap V+\op{Ric}(\grd V,\grd V)-V\lap V-|\grd V|^2 \\
&=|D^2V|^2+V^{-1}D^2V(\grd V,\grd V)-nV^2-|\grd V|^2 \\
&=|D^2V-Vg|^2-\frac{1}{2}V^{-1}\grd V\cdot\grd\(|\grd V|^2-V^2\).
\end{align*}
This can be rewritten as
\begin{equation}\label{dv-v}
\lap\(|\grd V|^2-V^2+1\)-V^{-1}\grd V\cdot\grd \(|\grd V|^2-V^2+1\)=2|D^2V+Vg|^2\geq 0.
\end{equation}

Suppose $(M,g,V)$ further satisfies conditions 1 and 2 in Theorem 1.
We consider the compact manifold with boundary $\bM$ with the Fermat
metric $\bg=(V+1)^{-2}g$ with.
\begin{lem}\label{confdef}
$(\bM,\bg)$ has scalar curvature $\bR\geq 0$. Moreover the boundary
is isometric to the standard sphere $S^{n-1}$ and has mean curvature $n-1$.
\end{lem}

\begin{prf}
By (\ref{bm}) we have
$$|\grd V|^2-V^2+1\ra 0, \text{\quad as} \ x\ra \S.$$
By the maximum principle for the equation (\ref{dv-v})
\begin{equation}\label{max}
|\grd V|^2-V^2+1\leq 0.
\end{equation}
The scalar curvature of $\bR$ of $\bg$ is given by
\begin{equation}
(V+1)^{-2}\bR=R+2(n-1)\frac{\lap V}{V+1}-n(n-1)\frac{|\grd V|^2}{(V+1)^2}
\end{equation}
By (\ref{max}) and the fact $R=-n(n-1)$
\begin{equation*}
(V+1)^{-2}\bR\geq -n(n-1)+2n(n-1)\frac{V}{V+1}-n(n-1)\frac{V^2-1}{(V+1)^2}
=0
\end{equation*}

By the asymptotic expansions (\ref{exph})(\ref{expv}) for $g$ and $V$, we
can write $\bg=\ph^{-2}(dr^2+h_r)$ near the boundary with
$\ph=r(V+1)=1+r+r^2/4+\ldots$. Then it's obvious that
$\bg|_{\S}=h_0$ is the round metric on $S^{n-1}$ and $\S$ has mean
curvature $\bH=n-1$.
\end{prf}

We now gives another proof of  Theorem 1 which does not require
the spin assumption for $n\leq 7$, \ie we have

\begin{thm}
Let $(M,g,V)$ be a static solution to the vacuum Einstein equation with
negative cosmological constant (\ref{va}). If
\begin{enumerate}
\item $(M,g)$ is conformally compact and the conformal boundary is $S^{n-1}$,
\item $V^{-1}$ is a defining function and $V^{-2}g|_{S^{n-1}}$ is the
      standard metric on $S^{n-1}$
\item either $n\leq 7$ or $M$ is spin,
\end{enumerate}
then $(M,g,V)$ is the anti-de Sitter spacetime.
\end{thm}

\begin{prf}
By Lemma \ref{confdef}, $(\bM,\bg=(1+V)^{-2}g)$ is a compact manifold with
boundary and the boundary is the standard sphere $S^{n-1}$ with mean
curvature $n-1$. Then $(\bM,\bg)$
is isometric to the unit ball $\bB\subset \R^n$. 
To see this, we attach $\R^n-\bB$ to $(\bM,\bg)$
along the boundary and get a complete manifold with a flat end.
Positive mass theorem for asymptotically manifolds
would imply that it is isometric to $\R^n$ if either $n\leq 7$ or $M$ is spin
except the metric is merely $C^1$ along the common boundary hypersurface.
It was known to many experts that the positive mass theorem holds
in this general setting and in fact this fact was used by Bunting 
and Masood-ul-alam \cite{BM, Mas1, Mas2} to prove uniqueness for various blackhole
solutions.
But there had been no detailed proofs in the literature until two recent papers by
Miao \cite{Miao} who uses a mollification argument to reduce it to the regular case and
by Shi-Tam \cite{ShT}  generalizing Witten's spinor argument.

Therefore $(\bM,\bg)$ is isometric
to $\bB$ and in particular $\bR=0$.
By the proof of Lemma \ref{confdef}, we have $D^2V=Vg$ and $|\grd V|^2-V^2+1=0$.
Then by (\ref{f1}) we have $\op{Ric}(g)=-(n-1)g$. As $g$ is also conformally
flat, it follows that $(M, g)$ is the hyperbolic space.
\end{prf}

\providecommand{\bysame}{\leavevmode\hbox to3em{\hrulefill}\thinspace}
\providecommand{\MR}{\relax\ifhmode\unskip\space\fi MR }
\providecommand{\MRhref}[2]{%
  \href{http://www.ams.org/mathscinet-getitem?mr=#1}{#2}
}
\providecommand{\href}[2]{#2}


\begin{thebibliography}{[99]}

\bibitem{BGH}
W.~Boucher, G.~W. Gibbons, and Gary~T. Horowitz, \emph{Uniqueness theorem for
  anti-de {S}itter spacetime}, Phys. Rev. D (3) \textbf{30} (1984), no.~12,
  2447--2451. \MR{86e:83014}

\bibitem{BM}
Gary~L. Bunting and A.~K.~M. Masood-ul Alam, \emph{Nonexistence of multiple
  black holes in asymptotically {E}uclidean static vacuum space-time}, Gen.
  Relativity Gravitation \textbf{19} (1987), no.~2, 147--154. \MR{88e:83031}

\bibitem{CH}
Piotr~T. Chru{\'s}ciel and Marc Herzlich, \emph{The mass of asymptotically
  hyperbolic riemannian manifolds}, e-Print math.DG/0110035.

\bibitem{ChS}
Piotr~T. Chru{\'s}ciel and Walter Simon, \emph{Towards the classification of
  static vacuum spacetimes with negative cosmological constant}, J. Math. Phys.
  \textbf{42} (2001), no.~4, 1779--1817. \MR{1 820 431}

\bibitem{Gra}
C.~Robin Graham, \emph{Volume and area renormalizations for conformally compact
  {E}instein metrics}, The Proceedings of the 19th Winter School ``Geometry and
  Physics'' (Srn\'\i, 1999), no.~63, 2000, pp.~31--42. \MR{2002c:53073}

\bibitem{GSW}
G.J. Galloway, S.~Surya, and E.~Woolgar, \emph{On the geometry and mass of
  static, asymptotically ads spacetimes, and the uniqueness of the ads
  soliton}, e-Print hep-th/0204081.


\bibitem{Miao}
Pengzi Miao, \emph{Positive Mass Theorem on Manifolds admitting Corners along
  a Hypersurface} arXiv math-ph/0212025.

\bibitem{Mas2}
A.~K.~M. Masood-ul Alam, \emph{Uniqueness proof of static charged black holes
  revisited}, Classical Quantum Gravity \textbf{9} (1992), no.~5, L53--L55.
  \MR{93f:83059}

\bibitem{Mas1}
\bysame, \emph{Uniqueness of a static charged dilaton black hole}, Classical
  Quantum Gravity \textbf{10} (1993), no.~12, 2649--2656. \MR{94k:83050}

\bibitem{Qing}
Jie Qing, \emph{On the rigidity for conformally compact {E}instein manifolds},
  Int. Math. Res. Not. (2003), no.~21, 1141--1153. \MR{1 962 123}

\bibitem{ShT}
Y. Shi and L. Tam, \emph{Positive mass theorem and the boundary behavior of 
compact manifolds with nonegative scalar curvature}, arXiv math.DG/0301047.

\bibitem{Wang2}
Xiaodong Wang, \emph{The mass of asymptotically hyperbolic manifolds}, J.
  Differential Geom. \textbf{57} (2001), no.~2, 273--299. \MR{2003c:53044}

\bibitem{Wang1}
\bysame, \emph{On conformally compact {E}instein manifolds}, Math. Res. Lett.
  \textbf{8} (2001), no.~5-6, 671--688. \MR{2003d:53075}

\end{thebibliography}
\end{document}